\documentclass[12pt,a4paper]{article}
\usepackage{}
\usepackage{indentfirst}
\usepackage{amsfonts}
\usepackage{amssymb}
\usepackage{mathrsfs}
\usepackage{amsmath}
\usepackage{amsthm}
\usepackage{enumerate}
\usepackage{cite}
\usepackage[all]{xy}
\allowdisplaybreaks
\newtheorem{defn}{Definition}[section]
\newtheorem{thm}[defn]{Theorem}
\newtheorem{lem}[defn]{Lemma}
\newtheorem{prop}[defn]{Proposition}
\newtheorem{cor}[defn]{Corollary}
\newtheorem{eg}[defn]{Example}
\newtheorem{re}[defn]{Remark}
\newcommand{\bdefn}{\begin{defn}}
\newcommand{\edefn}{\end{defn}}
\newcommand{\bthm}{\begin{thm}}
\newcommand{\ethm}{\end{thm}}
\newcommand{\blem}{\begin{lem}}
\newcommand{\elem}{\end{lem}}
\newcommand{\bprop}{\begin{prop}}
\newcommand{\eprop}{\end{prop}}
\newcommand{\bcor}{\begin{cor}}
\newcommand{\ecor}{\end{cor}}
\newcommand{\beg}{\begin{eg}}
\newcommand{\eeg}{\end{eg}}
\newcommand{\bre}{\begin{re}}
\newcommand{\ere}{\end{re}}
\newcommand{\bpf}{\begin{proof}}
\newcommand{\epf}{\end{proof}}

\newcommand{\supercite}[1]{\textsuperscript{\cite{#1}}}
\newcommand{\benu}{\begin{enumerate}}
\newcommand{\eenu}{\end{enumerate}}
\newcommand{\bc}{\begin{center}}
\newcommand{\ec}{\end{center}}
\newcommand{\bea}{\begin{eqnarray}}
\newcommand{\eea}{\end{eqnarray}}
\newcommand{\Bea}{\begin{eqnarray*}}
\newcommand{\Eea}{\end{eqnarray*}}
\newcommand{\beq}{\begin{equation}}
\newcommand{\eeq}{\end{equation}}
\newcommand{\Beq}{\begin{equation*}}
\newcommand{\Eeq}{\end{equation*}}
\newcommand{\bspl}{\begin{split}}
\newcommand{\espl}{\end{split}}
\newcommand\relphantom[1]{\mathrel{\phantom{#1}}}

\numberwithin{equation}{section}

\begin{document}
\title{\bf  Derivations from the even parts into the odd parts for Hamiltonian superalgebras}
\author{\normalsize \bf Yuan Chang, Liangyun Chen}
\date{\small{ School of Mathematics and Statistics, Northeast Normal University, Changchun,  130024,  CHINA}}
\maketitle
\begin{abstract}
Let $W_{\overline{1}}$ and $H_{\overline{0}}$ denote the odd parts of the general Witt modular Lie superalgebra $W$ and the even parts of the Hamiltonian Lie superalgebra $H$ over a field of characteristic $p>3$, respectively.
We give a torus of $H_{\overline{0}}$
and the weight space decomposition of the special subalgebra of $W_{\overline{1}}$ with respect to the torus.
By means of the derivations of the weight 0 and three series of outer derivations from $H_{\overline{0}}$ into $W_{\overline{1}}$, the derivations from the even parts of Hamiltonian superalgebra to the odd parts of Witt superalgebra are determined.

\bigskip
\noindent{Key words:}  Torus; Weight space decomposition; Derivations\\
\end{abstract}
\footnote[0]{Corresponding author(L. Chen): chenly640@nenu.edu.cn.}
\footnote[0]{Supported by  NNSF of China (Nos. 11171055 and No. 11471090). }
\section{Introduction}
As the natural generalization of Lie algebras, Lie superalgebras become an efficient tool for analyzing the properties of physical systems.
The theory of Lie superalgebras is closely related to many branches of mathematics.
In particular, V.G.Kac classified the finite-dimensional simple Lie superalgebras over algebraically closed field of characteristic zero\supercite{Kac}.
For modular Lie superalgebras, as we know, \cite{Kochetkov&Leites, Petrogradski} should be the earliest papers.
With the development of modular Lie superalgebras, Cartan-type Lie superalgebras play an important role in the category of modular Lie superalgebras.
Eight families of $\mathbb{Z}$-graded Cartan-type Lie superalgebras were constructed over a field of characteristic $p>3$\supercite{Zhang, Fu&Zhang&Jing, Liu&He, Liu&Yuan}.
Determining the superderivation algebras are very important subjects in modular Lie superalgebras.
The superderivation algebras were studied in one-by-one fashion for the finite dimensional and
simple ones\supercite{Zhang, Fu&Zhang&Jing, Liu&He, Liu&Yuan, Ma&Zhang, Wang&Zhang, Zhang&Zhang}.
\cite{Bai&Liu} used a uniform method to determine the superderivation algebras of $\mathbb{Z}$-graded Cartan-type Lie superalgebras.

For Lie superalgebras, the even parts are closely connected with Lie algebras and the odd parts are the modulars of the even parts.
Determining the derivations of the even parts and the derivations from the even parts into the odd parts for Lie superalgebras are very interesting.
\cite{Liu&Zhang, Liu&Guan} respectively determine the derivations of the even parts and the derivations from the even parts into the odd parts for Lie superalgebras $W$ and $S$ of Cartan type.
\cite{Guan&Liu, Guan&Chen} respectively determine the derivations of the even parts and the derivations from the even parts into the odd parts for contact Lie superalgebras.
The odd $\mathbb{Z}$-homogeneous derivations and negative $\mathbb{Z}$-homogeneous derivations from the even parts of Hamiltonian Lie superalgebras into the even parts of Witt Lie superalgebras\supercite{Liu&Su&Zhang}.
Moreover, there are more outer derivations for the even parts of Hamiltonian Lie superalgebras than for the even parts of Lie superalgebras $W$ and $S$.
The purpose of this paper is to determine the derivations from the even parts of the Hamiltonian Lie superalgebras into the odd parts $W_{\overline{1}}$ of the Witt Lie superalgebras.

This paper is organized as follow.
In Section 2, we give the basic notations and concepts.
In Section 3, we introduce 
the suitable generating set and the proper torus $T$ of the even parts $\mathcal{H}$ in the finite-dimensional Hamiltonian Lie superalgebras $H$.
Then we give the special subalgebra $\mathfrak{g}_1$ of $W_{\overline{1}}$ and the weight space decompositions of $\mathfrak{g}_1$ with respect of the torus $T$ that have the same weights with some generators of $\mathcal{H}$.
In Section 4, we give three series of outer derivations from the even parts $H_{\overline{0}}$ into the odd parts $W_{\overline{1}}$ .
We characterize the derivations vanishing on the top of $H_{\overline{0}}$ with the above results.
In Section 5, we determine the derivation algebras from the even parts of Hamiltonian Lie superalgebras to the odd parts of Witt modular Lie superalgebras.

\section{Basic}
For a vector superspace $V=V_{\overline{0}}\oplus V_{\overline{1}}$,~we denote by $\mathrm{p}(x)=\alpha$ the parity of a homogeneous element $x\in V_{\alpha}$,~$\alpha\in \mathbb{Z}_2$, where $\mathbb{Z}_2=\{\bar{0},\bar{1}\}$.~
If $V=\oplus_{i\in\mathbb{Z}} V_i$ is a $\mathbb{Z}$-graded vector space and $x\in V$ is a $\mathbb{Z}$-homogeneous element,~write $\mathrm{zd}(x)$ for the $\mathbb{Z}$-degree of $x$.~The symbol $\mathrm{p}(x)$ (resp.$\mathrm{zd}(x)$) always implies that $x$ is a
$\mathbb{Z}_2$-(resp.$\mathbb{Z}$-)homogeneous element.~
If $L=\oplus_{-r\leq i\leq s} L_i$ is a $\mathbb{Z}$-graded Lie superalgebra,
then $\oplus_{-r\leq i\leq0} L_i$ is called the $top$ of $L$.


Now,~we review the notions of modular Lie superalgebras $W$ and $H$ of Cartan-type and their grading structures.
Throughout $\mathbb{F}$ is an algebraically closed field of characteristic $p>3$.
We write $\mathbb{N}$ for the set of natural numbers and~$\mathbb{N}_{+}$ for the set of positive integers.
Fix two positive integers $m$ and $n$.
For $\alpha=(\alpha_1,\cdots,\alpha_m)\in\mathbb{N}_{+}^{m}$,~put $|\alpha|=\Sigma_{i=1}^{m}\alpha_i$.
In \cite{Strade}, denote by $\mathcal{O}(m)$ the $divided~power~algebra$ over $\mathbb{F}$ with an $\mathbb{F}$-basis $\{x^{(\alpha)}|\alpha\in\mathbb{N}_{+}^{m}\}$.
For $\varepsilon_i=(\delta_{i1},\cdots,\delta_{im})$,~we abbreviate $x^{(\varepsilon_i)}$ to $x_i$,~$i=1,\cdots,m$.
Let $\Lambda(n)$ be the $exterior~superalgebra$ over $\mathbb{F}$ in $n$ variables $x_{m+1},\cdots,x_{m+n}$.
Denote the tensor product by $\mathcal{O}(m,n)=\mathcal{O}(m)\otimes_{\mathbb{F}}\Lambda(n)$.
Obviously,~$\mathcal{O}(m,n)$ is an associative superalgebra with a $\mathbb{Z}_2$-gradation induced by the trivial $\mathbb{Z}_2$-gradation of $\mathcal{O}(m)$ and the natural $\mathbb{Z}_2$-gradation of $\Lambda(n)$.
Moreover,~$\mathcal{O}(m,n)$ is super-commutative.
For $g\in\mathcal{O}(m)$,~$f\in\Lambda(n)$,~we write $gf$ for $g\otimes f$.
The following formulas hold in $\mathcal{O}(m,n)$:
$$x^{(\alpha)}x^{(\beta)}={\alpha+\beta\choose\alpha}x^{(\alpha+\beta)}~~~~~~~~~~~~~\mathrm{for}~\alpha,\beta\in\mathbb{N}_{+}^{m};$$
$$x_kx_l=-x_lx_k~~~~~~~~~~~\mathrm{for}~k,l=m+1,\cdots,m+n;$$
$$x^{(\alpha)}x_k=x_kx^{(\alpha)}~~~~~\mathrm{for}~\alpha\in\mathbb{N}_0^{m}, k=m+1,\cdots,m+n,$$
where ${\alpha+\beta\choose\alpha}:=\Pi_{i=1}^{m}{\alpha_i+\beta_i\choose\alpha_i}$.
Put~$\mathrm{I}_0:=\{1,2,\cdots,m\}$,~$\mathrm{I}_1:=\{m+1,\cdots,m+n\}$~and~$\mathrm{I}:=\mathrm{I}_0\cup \mathrm{I}_1$.
Set
$$\mathbb{B}_k:=\{\langle i_1,i_2,\cdots,i_k\rangle\mid m+1\leq i_1<i_2<\cdots<i_k\leq m+n\}$$
and~$\mathbb{B}:=\mathbb{B}(n)=\cup_{k=0}^{n}\mathbb{B}_k$,~where~$\mathbb{B}_0=\emptyset$.
For~$u=\langle i_1,i_2,\cdots,i_k\rangle\in\mathbb{B}_k$,~set~$|u|:=k$,~$x^u=x_{i_1}\cdots x_{i_k}$,~$|\emptyset|=0$ and $x^{\emptyset}=1$.
For convenience,~we use~$\langle s\rangle$ to stand for~$\langle\delta_{1s}, \delta_{2s},\cdots, \delta_{(m+n)s}\rangle$, for $s\in\mathrm{I}_1 $.
Set~$\mathbb{B}^0=\{u\in\mathbb{B}||u|$ is even$\}$. 
Clearly,~$\{x^{(\alpha)}x^u|\alpha\in\mathbb{N}_{+}^{m},u\in\mathbb{B}\}$~constitutes an $\mathbb{F}$-basis of $\mathcal{O}(m,n)$.
Let~$\partial_1,\partial_2,\cdots,\partial_{m+n}$ be the linear transformations of~$\mathcal{O}(m,n)$ such that $\partial_i(x^{(\alpha)})=x^{(\alpha-\varepsilon_i)}$ for $i\in \mathrm{I}_0$
and $\partial_i(x_k)=\delta_{ik}$ for $i\in \mathrm{I}_1$.
Obviously,~$\mathrm{p}(\partial_i)=\mu(i)$,~where~$\mu(i):=\overline{0}$,~$i\in \mathrm{I}_0$,~and~$\mu(i):=\overline{1}$,~$i\in \mathrm{I}_1$.
Then $\partial_1,\partial_2,\cdots\partial_{m+n}$ are superderivations of the superalgebra $\mathcal{O}(m,n)$.
Let
$$W(m,n):=\left\{\sum_{r\in I}f_r\partial_r|f_r\in\mathcal{O}(m,n),r\in I\right\}.$$
Then $W(m,n)$ is an infinite-dimensional Lie superalgebra contained in $\mathrm{Der}(\mathcal{O}(m,n))$.
One can verify that $$[fD,gE]=fD(g)E-(-1)^{\mathrm{p}(fD)\mathrm{p}(gE)}gE(f)D+(-1)^{\mathrm{p}(D)\mathrm{p}(g)}fg[D,E]$$
for~$f,g\in\mathcal{O}(m,n)$,~$D,E\in\mathrm{Der}(\mathcal{O}(m,n))$.
Specially,
$$[f\partial_i,g\partial_j]=f\partial_i(g)\partial_j-(-1)^{\mathrm{p}(f\partial_i)\mathrm{p}(g\partial_j)}g\partial_j(f)\partial_i~~~~\mathrm{for} ~f,~g\in\mathcal{O}(m,n),~ i,j\in \mathrm{I}.$$

Hereafter,~suppose~that $m=2r,~r\in\mathbb{N}_0$.
If $n$ is even, then we set~$n=2s$,~otherwise~$n=2s+1$, where $s\in\mathbb{N}_0$.
Define a linear mapping $\mathrm{D_H}:\mathcal{O}(m,n)\rightarrow W(m,n)$~by means of
$$\mathrm{D_H}(f):=\sum_{i\in I}\tau(i)(-1)^{\mu(i)\mathrm{p}(f)}\partial_i(f)\partial_{i'}~~~~~~~~\mathrm{for~all}~f\in\mathcal{O}(m,n),$$
where
\begin{equation*}
i'=
\begin{cases}
i+r,&~~~~~1\leq i\leq r,
\\i-r,&~~~~~r+1\leq i\leq m,
\\i+s,&~~~~~m+1\leq i\leq m+s,
\\i-s,&~~~~~m+s<i\leq m+2s,
\\i,&~~~~~the~other,
\end{cases}~~~~~~~
\tau(i)=
\begin{cases}
1,&~~1\leq i\leq r,
\\-1,&~~r<i\leq m,
\\1,&~~i\in\mathrm{I}_1.
\end{cases}
\end{equation*}
The following equation holds:$$[\mathrm{D_H}(f),\mathrm{D_H}(g)]=\mathrm{D_H}(\mathrm{D_H}(f)(g))~~~~~~~~\mathrm{for~all}~f,~g\in\mathcal{O}(m,n).$$
Put $H(m,n):=\mathrm{span}_{\mathbb{F}}\{\mathrm{D_H}(f)|f\in\mathcal{O}(m,n)\}$.
Obviously,~$H(m,n)$ is an infinite-dimensional $\mathbb{Z}_2$-graded subalgebra of $W(m,n)$.

Fix two $m$-tuples of positive integers: $$\underline{t}:=(t_1,t_2,\cdots,t_m)\in\mathbb{N}^m,~~~~~~\pi:=(\pi_{1},\pi_{2},\cdots,\pi_{m}),$$
where $\pi_i=p^{t_{i}}-1$,~$i=1,2,\cdots,m$.~
Let $\mathbb{A}:=\mathbb{A}(m;\underline{t})=\{\alpha\in\mathbb{N}_0^{m}|\alpha_i\leq\pi_i,~i\in \mathrm{I}_0\}$.~
Then $$\mathcal{O}(m,n;\underline{t}):=\mathrm{span}_{\mathbb{F}}\left\{x^{(\alpha)}x^u|~\alpha\in\mathbb{A},~u\in\mathbb{B}\right\}$$
is a finite-dimensional subalgebra of~$\mathcal{O}(m,n)$~with a $\mathbb{Z}$-grading structure: $$\mathcal{O}(m,n;\underline{t})=\bigoplus_{i=0}^{\xi}\mathcal{O}(m,n;\underline{t})_i,$$
where $\mathcal{O}(m,n;\underline{t})_i:=\mathrm{span}_{\mathbb{F}}\{x^{(\alpha)}x^u|~|\alpha|+|u|=i\}~\mathrm{and}~\xi:=|\pi|+n .$
Set
$$W(m,n;\underline{t}):=\left\{\sum_{i\in \mathrm{I}}f_i\partial_i|f_i\in\mathcal{O}(m,n;\underline{t}),~i\in \mathrm{I}\right\}.$$
Then~$W(m,n;\underline{t})$ is a finite-dimensional simple Lie superalgebra\supercite{Zhang},~
which is called the generalized Witt Lie superalgebras.~
We note that $W(m,n;\underline{t})$ possesses a standard $\mathbb{F}$-basis $\{x^{(\alpha)}x^u\partial_r|\alpha\in\mathbb{A},u\in\mathbb{B},r\in \mathrm{I}\}$
and a $\mathbb{Z}$-grading structure:
$$W(m,n;\underline{t})=\bigoplus_{r=-1}^{\xi-1}W(m,n;\underline{t})_r,$$
where $W(m,n;\underline{t})_r:=\mathrm{span}_{\mathbb{F}}\{x^{(\alpha)}x^u\partial_j|~|\alpha|+|u|=r+1,~j\in\mathrm{I}\}.$
Set
$$H(m,n;\underline{t}):=\left\{\mathrm{D_H}(f)|f\in\bigoplus_{i=0}^{\xi-1}\mathcal{O}(m,n;\underline{t})_i\right\},$$
where $\xi:=|\pi|+n$.
Then~$H(m,n;\underline{t})$ is a finite-dimensional $\mathbb{Z}$-graded simple subalgebra of $W(m,n;\underline{t})$\supercite{Zhang},~which is called the Hamiltonian superalgebra.~

In the following sections,~$\mathcal{O}(m,n;\underline{t})$,~$W(m,n;\underline{t})$~and~$H(m,n;\underline{t})$ will be denoted by $\mathcal{O}$,~$W$~and $H$.
In addition,~the even parts of $W$ and $H$ will respectively be denoted by $W_{\overline{0}}$ and $H_{\overline{0}}$,~and the odd parts of $W$ 
will be denoted by $W_{\overline{1}}$. 

\section{Torus and Weight Space Decompositions}
Recall 
$$H_{\overline{0}}=\mathrm{span}_{\mathbb{F}}\{\mathrm{D_H}(x^{(\alpha)}x^u)| \alpha\in\mathbb{A},u\in\mathbb{B}^0,(\alpha,u)\neq(\pi,\omega)\}.$$
Set$$\mathcal{J}:=\mathrm{span}_{\mathbb{F}}\{\mathrm{D_H}(x^{(\alpha)}x^u)|\alpha\in\mathbb{A},u\in\mathbb{B}^0,(\alpha,u)\neq(\pi,\omega),(\alpha,u)\neq(\pi,\emptyset)\},$$
where $\mathbb{B}^0=\{u\in\mathbb{B}||u|$ is even$\}$ and $\omega=\langle m+1,m+2,\cdots,m+n\rangle\in\mathbb{B}_n$.~
Obviously,~$\mathcal{J}$~is a subspace of $H_{\overline{0}}$ of codimension 1:
$$H_{\overline{0}}=\mathcal{J}\oplus\mathbb{F}\mathrm{D_H}(x^{(\pi)}).$$

\blem$^{[16]}$
$\mathcal{J}$ is a maximal idea of $H_{\overline{0}}$.
\elem

We review the generating sets of $\mathcal{J}$.
Set $$\mathcal{M}=\{\mathrm{D_H}(x^{(q_i\varepsilon_i)})|1\leq q_i\leq\pi_i,~i\in\mathrm{I}_0\}~~\mathrm{and}~~\mathcal{N}=\{\mathrm{D_H}(x_ix^u)|i\in\mathrm{I}_0,~u\in\mathbb{B}_2\}.$$

\blem$^{[16]}$
$\mathcal{J}$ is generated by $\mathcal{M}\cup\mathcal{N}\cup\mathcal{J}_0$.
\elem
Set $\mathfrak{g}_1=\mathrm{span}_{\mathbb{F}}\{x^u\partial_r|r\in\mathrm{I},~u\in\mathbb{B},~\mathrm{p}(x_u\partial_r)=\overline{1}\}$.
And $\mathfrak{g}_1$ is a $\mathbb{Z}$-graded subspace of $W_{\overline{1}}$.
Since $C_{W_{\overline{1}}}((W_{\overline{0}})_{-1})=\mathfrak{g}_1$,
$C_{W_{\overline{1}}}(\mathcal{J})\subseteq C_{W_{\overline{1}}}((H_{\overline{0}})_{-1})=C_{W_{\overline{1}}}((W_{\overline{0}})_{-1})=\mathfrak{g}_1$.

\bprop\label{lem3}
The following statements hold:
\begin{itemize}
\item[\rm{(1)}] If $|\omega|=n$ is even,~then $C_{W_{\overline{1}}}(\mathcal{J})=0$.
\item[\rm{(2)}] If $|\omega|=n$ is odd,~then $C_{W_{\overline{1}}}(\mathcal{J})=\mathbb{F}\mathrm{D_H}(x^{\omega})$.
\end{itemize}
\eprop

\bpf
Since~$C_{W_{\overline{1}}}(\mathcal{J})$ is a $\mathbb{Z}$-graded subalgebra of $W_{\overline{1}}$,~
we note that $C_{W_{\overline{1}}}(\mathcal{J})\subseteq\mathfrak{g}_1$.~
If $D\in C_{W_{\overline{1}}}(\mathcal{J})$, we see that $D\in\mathfrak{g}_1$.~
Thus one may assume that
\begin{eqnarray}\label{eq1}
D=\sum_{r\in\mathrm{I}}f_r\partial_r,
\end{eqnarray}
where~$f_r\in\Lambda(n)$.~
For any $i\in\mathrm{I}_0$,~since~$\mathrm{D_H}(x^{(2\varepsilon_i)})\in\mathcal{J}$,~
we have $[D,\mathrm{D_H}(x^{(2\varepsilon_i)})]=f_i\partial_i$ and therefore $f_i\partial_i=0$.~
This prove that $f_i=0$ for all $i\in\mathrm{I}_0$.~
By Eq.(\ref{eq1}), we can obtain that
\begin{eqnarray}\label{eq2}
D=\sum_{r\in\mathrm{I}_1}f_r\partial_r,~~~~~where~~f_r\in\Lambda(n).
\end{eqnarray}

%
%
$\rm{(1)}$~~When $|\omega|=n$ is even,~there are three cases to discuss.~

Case 1. If $2\leq\mathrm{zd}(f_r)\leq n-2$~for $r\in\mathrm{I}_1$,
assume that $f_{k}\neq0$ for some $k\in\mathrm{I}_1$.
We choose $k,~k'\in\mathrm{I}_1$.~
Then
\begin{eqnarray}\label{eq3}
0=[\mathrm{D_H}(x^kx^{k'}),D]=\sum_{r\in\mathrm{I}_1}(x^{k'}\partial_{k'}(f_r)-x^{k}\partial_{k}(f_r))\partial_r-f_{k'}\partial_{k'}+f_{k}\partial_{k}.
\end{eqnarray}
By comparison the coefficients of $\partial_k$ in Eq.(\ref{eq3}),~we can find $x^k\in f_k$ and
$x^{k'}\notin f_k$.
For $l,~l'\in\mathrm{I}_1\backslash\{k,~k'\}$,~we have
\begin{eqnarray}\label{eq4}
0=[\mathrm{D_H}(x^lx^{l'}),D]=\sum_{r\in\mathrm{I}_1}(x^{l'}\partial_{l'}(f_r)-x^{l}\partial_{l}(f_r))\partial_r-f_{l'}\partial_{l'}+f_{l}\partial_{l}.
\end{eqnarray}
Since the coefficients of $\partial_k$ in Eq.(\ref{eq4}) yield that $x^{l'}\partial_{l'}(f_k)-x^l\partial_{l}(f_k)=0$,
we can obtain that $x^l,x^{l'}\in f_k$~or $x^l,x^{l'}\notin f_k$.
So $\mathrm{p}(f_k)$ is odd and $\mathrm{p}(f_k\partial_k)$ is even,
which is a contradiction.
This proves that $D=0$.

Case 2. If $\mathrm{zd}(f_r)=0~\mathrm{for}~r\in\mathrm{I}_1$,~
then we assume that
$$D=\sum_{r\in\mathrm{I}_1}a_r\partial_r,~~~~~~~~\mathrm{where}~~a_r\in\mathbb{F}.$$
For any $k,~l\in\mathrm{I}_1$,~we obtain
\begin{eqnarray}\label{eq5}
0=[\mathrm{D_H}(x^kx^l),\sum_{r\in\mathrm{I}_1}a_r\partial_r]=-a_l\partial_{k'}+a_k\partial_{l'}.
\end{eqnarray}
By comparison the coefficients in Eq.(\ref{eq5}),~we obtain that $a_r=0$ for any $r\in\mathrm{I}_1$.
It is obvious that $D=0$.

Case 3. If $\mathrm{zd}(f_r)=n~\mathrm{for}~r\in\mathrm{I}_1$.~
We can assume
$$D=\sum_{r\in\mathrm{I}_1}a_{r}x^{\omega}\partial_r,~~~~~~\mathrm{where}~~a_r\in\mathbb{F}.$$
Arguing as the proof of Case 2,~we can obtain $a_r=0$ for any $r\in\mathrm{I}_1$ and $D=0$.

$\rm{(2)}$~~
We declare $|w|=n$ is odd.~
For any basis element $\mathrm{D_H}(x^{(\alpha)}x^u)$ of $\mathcal{J}$,
it is clear that $[\mathrm{D_H}(x^{\omega}),\mathrm{D_H}(x^{(\alpha)}x^u)]=0$.
Therefore,~$\mathbb{F}\mathrm{D_H}(x^{\omega})\subseteq C_{W_{\overline{1}}}(\mathcal{J})$.
We will show the converse inclusion.
Every element $D\in C_{W_{\overline{1}}}(\mathcal{J})$ may be written as Eq.(\ref{eq2}).
We distinguish two cases to discuss.

Case 1. $0\leq\mathrm{zd}(f_r)\leq n-3~\mathrm{for}~r\in\mathrm{I}_1$.
Arguing just as the proof of the above Case 1 and Case 2,~we can proof that $D=0$.

Case 2. $\mathrm{zd}(f_r)=n-1~\mathrm{for}~r\in\mathrm{I_1}$.
Then one may assume that
$$f_r=\sum_{s\in\mathrm{I_1}}c_{rs}x^{\widehat{s}},~~~~~~~~~~\mathrm{where}~c_{rs}\in\mathbb{F}.$$
Then we put~$\widehat{s}:=\omega- \langle s\rangle$ where ~for $s\in\mathrm{I_1}$.~
For $k,~k'\in\mathrm{I}_1$,~we have the equation that
$0=[\mathrm{D_H}(x^kx^{k'}),D]$.
Then
$$0=[x^{k'}\partial_{k'}-x^{k}\partial_{k},\sum_{r\in\mathrm{I_1}}(\sum_{s\in\mathrm{I_1}}c_{rs}x^{\widehat{s}})\partial_r]=\sum_{r\in\mathrm{I_1}}(c_{rk}x^{\widehat{k}}-c_{rk'}x^{\widehat{k'}})\partial_r
+(\sum_{s\in\mathrm{I_1}}c_{ks}x^{\widehat{s}})\partial_k
-(\sum_{s\in\mathrm{I_1}}c_{k's}x^{\widehat{s}})\partial_{k'}.$$
By the calculation of the coefficients of $\partial_r$ for $r\in \mathrm{I_1}$,
we obtain that 
$$D=\sum_{r\in\mathrm{I_1}}c_{rr'}x^{\widehat{r'}}\partial_r.$$
For any $k,~l\in\mathrm{I_1}$,~we have the equation that
$$0=[\mathrm{D_H}(x^kx^{l}),D]=[x^{l}\partial_{k'}-x^{k}\partial_{l'},\sum_{r\in\mathrm{I_1}}c_{rr'}x^{\widehat{r'}}\partial_r].$$
We obtain that $(c_{kk'}+(-1)^{k'+l-3}c_{l'l})x^{\widehat{k'}}\partial_{l'}-(c_{ll'}+(-1)^{l'+k-3}c_{k'k})x^{\widehat{l'}}\partial_{k'}=0$,~
and therefore, $c_{k'k}=(-1)^{k+l}c_{l'l}$ for $k,~l\in\mathrm{I_1}$.~
Set $\mu=c_{i'i}$ where $i=m+1$.~
So far, we prove that
$$D=\mu\sum_{r\in\mathrm{I_1}}(-1)^{r+(m+1)}x^{\widehat{r}}\partial_{r'}=\mu\mathrm{D_{H}}(x^{\omega}).$$
\epf

%
%
\begin{thm}\label{lem4}
Suppose $\phi\in\mathrm{Der}(H_{\overline{0}},W_{\overline{1}})$ and $\phi(\mathcal{J})=0$.
Then the following statements hold:
\begin{itemize}
\item[\rm{(1)}] If $|\omega|=n$ is even,~then $\phi=0$.
\item[\rm{(2)}] If $|\omega|=n$ is odd,~then $\phi(\mathrm{D_{H}}(x^{(\pi)}))=\mu\mathrm{D_{H}}(x^{\omega})$
for some $\mu\in\mathbb{F}$.~Conversely,~any linear mapping $\phi:H_{\overline{0}}\rightarrow W_{\overline{1}}$~vanishing on $\mathcal{J}$ and satisfying $\phi(\mathrm{D_{H}}(x^{(\pi)}))=\mu(\mathrm{D_{H}}(x^{\omega}))$ for any fixed $\mu\in\mathbb{F}$ is a derivation from $H_{\overline{0}}$ into~$W_{\overline{1}}$.
\end{itemize}
\end{thm}

\bpf
By Lemma 3.1 and $[\phi(\mathrm{D_H}(x^{(\pi)})),\mathcal{J}]=0$,
that is $\phi(\mathrm{D_H}(x^{(\pi)}))\in C_{W_{\overline{1}}}(\mathcal{J})$,
we can get the consequences directly.
\epf

For any fixed $\mu\in\mathbb{F}$, define the linear mapping $\Gamma_{\mu}:H_{\overline{0}}\rightarrow W_{\overline{1}}$ by means of $\Gamma_{\mu}(\mathcal{J})=0$ and $ \Gamma_{\mu}(\mathrm{D_{H}}(x^{(\pi)}))=\mu\mathrm{D_{H}}(x^{\omega})$.
In the case that $|\omega|=n$ is odd, by Theorem \ref{lem4}, $\Gamma_{\mu}$ are outer derivations from $H_{\overline{0}}$ into~$W_{\overline{1}}$,~where $\mu\in\mathbb{F}$,
and $\mathrm{zd}(\Gamma_{\mu})=|\omega|-|\pi|$.

Suppose $L$ is a Lie superalgebra and $V$ an $L$-module.
Denote by $\mathrm{Der}(L,V)$ the $superderivation~space$ and $\mathrm{Ind}(L,V)$ the $inner~derivation~ space$.
Clearly, $\mathrm{Der}(L,V)$ is an $L$-submodule of $\mathrm{Hom}_{\mathbb{F}}(L,V)$.
A derivation $D:L\rightarrow V$ is called $inner$ if there is $v\in V$ such that $D(x)=xv$ for all $x\in L$.
Assume in addition that $L=\oplus_{r\in\mathbb{Z}}L_r$ is $\mathbb{Z}$-graded and finite-dimensional,
and $V=\oplus_{r\in\mathbb{Z}}V_r$ is a $\mathbb{Z}$-graded $L$-module.
Then the superderivation space inherits a $\mathbb{Z}$-graded $L$-module structure
$$\mathrm{Der}(L,V)=\bigoplus_{r\in\mathbb{Z}}\mathrm{Der}_r(L,V).$$

Let $(L,[p])$ be a restricted Lie algebra.
An element $x\in L$ is $p$-semisimple provided that $x\in \Sigma_{r\in\mathbb{N}}\mathbb{F}x^{[p]^r}$.
An abelian restricted subalgebra $T$ of $L$ is called a $torus$ if every element in $T$ is $p$-semisimple.
Let $T\subseteq L_0\cap L_{\overline{0}}$ be a torus of $L$ with the weight space decomposition:
$$L=\bigoplus_{\alpha\in\Theta}L_{\alpha},~~~~~~~V=\bigoplus_{\beta\in\Delta}V_{\beta}.$$
Then there exist subsets $\Theta_i\subset\Theta$~and $\Delta_j\subset\Delta$~such that
$L_i=\oplus_{\alpha\in\Theta_i}L_i\cap L_{\alpha}$~and $V_j=\oplus_{\beta\in\Delta_j}V_j\cap V_{\beta}$.
Hence $L$ and $V$ have the corresponding $\mathbb{Z}\times T^{\ast}$-grading structures, respectively,
where $T^{\ast}$ is the dual space of $T$.
Of course $\mathrm{Der}(L, V)$~inherits a $\mathbb{Z}\times T^{\ast}$-grading from $L$ and $V$ as above.
A superderivation $\phi\in\mathrm{Der}(L, V)$ is called a $weight$-$derivation$ if it is $T^{\ast}$-homogeneous.
Every superderivation is a sum of weight-derivations.

Set $T=\mathrm{span}_{\mathbb{F}}\{\mathrm{D_{H}}(x_ix_{i'})|i\in(1,\cdots,r)\cup(m+1,\cdots,m+s)\}$.
Obviously,~$T$ is a torus of $H_{\overline{0}}$.~
For any $\mathrm{D_H}(x_ix_{i'})\in T$,~we have
$$\mathrm{D_H}(x_ix_{i'})^p=\mathrm{D_H}(x_ix_{i'}),$$
\begin{eqnarray}\label{eq6}
[\mathrm{D_H}(x_ix_{i'}),\mathrm{D_H}(x^{(\alpha)}x^{u})]=(\alpha_{i'}\delta_{i'\in\mathrm{I_0}}-\alpha_{i}\delta_{i\in\mathrm{I_0}}+\delta_{i'\in u}-\delta_{i\in u})\mathrm{D_H}(x^{(\alpha)}x^{u}).
\end{eqnarray}
For $\alpha\in\mathbb{A}$ and $u\in \mathbb{B}$,~define a linear function $(\alpha+u)$ on $T$ such that
$$(\alpha+u)(\mathrm{D_H}(x_ix_{i'}))=(\alpha_{i'}\delta_{i'\in\mathrm{I_0}}-\alpha_{i}\delta_{i\in\mathrm{I_0}}+\delta_{i'\in u}-\delta_{i\in u}).$$
Further,~$H_{\overline{0}}$~and $W_{\overline{1}}$ both have weight space decompositions about $T$:
$$H_{\overline{0}}=\bigoplus_{(\alpha+u)}H_{\overline{0}~(\alpha+u)},~~~~~~~~W_{\overline{1}}=\bigoplus_{(\alpha+u)}W_{\overline{1}~(\alpha+u)}.$$

\blem$^{[13]}$\label{lem5}
Suppose that $L$ is a $\mathbb{Z}$-graded subalgebra of~$W_{\overline{0}}$ and $L_{-1}=(W_{\overline{0}})_{-1}$.
Let $E\in L$ and $\phi\in \mathrm{Der}(L,W_{\overline{1}})$ such that $\phi((W_{\overline{0}})_{-1})=0$.
Then $\phi(E)\in\mathfrak{g}_1$ if and only if $[E,(W_{\overline{0}})_{-1}]\subseteq\mathrm{ker}\phi$.
\elem

\blem$^{[11]}$\label{lem6}
A weight-derivation $\phi\in\mathrm{Der}(L,V)$ is inner if it is a nonzero weight-derivation.
In particular,~any derivation $\phi\in\mathrm{Der}(L,V)$ is inner modulo a zero weight-derivation.
\elem

\blem\label{lem7}
Let~$i,~j\in \mathrm{I_0}$~and $k,~l\in\mathrm{I_1}$.
For~$q_{i}^1$,~$q_{i}^2\in \mathbb{F}$,~set~$q_i^1\equiv0~(\mathrm{mod}p)$,~$%
q_i^2\equiv1~\linebreak(\mathrm{mod}p)$.~
Then the following statements hold:
\begin{itemize}
\item[\rm{(1)}] $\mathfrak{g}_{1(q_{i}^1\varepsilon_i)}=\left(\frac{1+(-1)^{n+1}}{2}\right)\left(\sum_{r\in\mathrm{I_1}}                   \mathbb{F}x^{\widehat{r'}}\partial_r\right)$.
\item[\rm{(2)}] $\mathfrak{g}_{1(q_{i}^2\varepsilon_i)}=\left(\frac{1+(-1)^{n+1}}{2}\right)\left( \mathbb{F}x^{\overline{u}}x^{m+n}\partial_{i'}\right)$.
\item[\rm{(3)}] $\mathfrak{g}_{1(\varepsilon_i+\langle k,l\rangle)}=\left(\frac{1+(-1)^{n+1}}{2}\right)\left( \mathbb{F}x^{k}x^{l}x^{m+n}\partial_{i'}\right)$.
\item[\rm{(4)}] $\mathfrak{g}_{1(\varepsilon_i\varepsilon_{i'})}=\left(\frac{1+(-1)^{n+1}}{2}\right)\left( \sum_{r\in\mathrm{I_1}}\mathbb{F}x^{\widehat{r'}}\partial_r\right)$.
\item[\rm{(5)}] $\mathfrak{g}_{1(\varepsilon_i\varepsilon_{j})}=\left(\frac{-1+(-1)^n}{2}\right)\left( \mathbb{F}x^{\overline{u}_{i'}}x^{m+n}\partial_{i'}+\mathbb{F}x^{\overline{u}_{j'}}x^{m+n}\partial_{j'}\right)$.
\item[\rm{(6)}] $\mathfrak{g}_{1(\varepsilon_k\varepsilon_{k'})}=\left(\frac{1+(-1)^{n+1}}{2}\right)\left( \sum_{r\in\mathrm{I_1}}\mathbb{F}x^{\widehat{r'}}\partial_r\right).$
\item[\rm{(7)}] $\mathfrak{g}_{1(\langle k,l\rangle)}=
    \mathbb{F}x^{\overline{u}_{k'}}\partial_{k'}+\mathbb{F}x^{\overline{u}_{l'}}\partial_{l'}.$
\end{itemize}
Where $x^l,~x^{l'}$ are both in $x^{\overline{u}}$ for $l\in\{m+1,\cdots,m+s\}$.
\elem

\bpf
Compare Eq.(\ref{eq6}) with the following equation for $i\in(1,\cdots,r)\cup(m+1,\cdots,m+s)$,
$$[\mathrm{D_H}(x_ix_{i'}),x^u\partial_r]=\delta_{i'\in u}-\delta_{i\in u}-\delta_{i'r}+\delta_{ir},~~~~~~~~\mathrm{where}~x^u\partial_r\in \mathfrak{g}_1.$$
By calculation, we can obtain the equations (1)$-$(7).
\epf

\section{Derivation Algebras}

In this section,~we determine the derivations from $\mathcal{J}$ into $W_{\overline{1}}$ which vanish on the top of $\mathcal{J}$.
To this aim, one needs to investigate the action of the derivations on the generators of $\mathcal{J}$.
We shall use the fact that $\mathcal{J}$ is generated by $\mathcal{M}\cup\mathcal{N}\cup\mathcal{J}_0$ (see Lemma 3.2).
Recall that
$$\mathfrak{g}_1=C_{W_{\overline{1}}}((W_{\overline{0}})_{-1})=\mathrm{span}_{\mathbb{F}}\{x^u\partial_r|r\in\mathrm{I},~u\in\mathbb{B},~\mathrm{p}(x_u\partial_r)=\overline{1}\}.$$
For simplicity, put
$$\mathrm{E}(\mathfrak{g}_1)=\bigoplus_{r\in\mathbb{Z}}\mathfrak{g}_{1(2r)},~~~~~~~~~ \mathrm{O}(\mathfrak{g}_1)=\bigoplus_{r\in\mathbb{Z}}\mathfrak{g}_{1(2r+1)},$$
where $\mathfrak{g}_{1(2r)}=\{x^u\partial_i||u|=2r+1, i\in\mathrm{I}_0\}$ and $\mathfrak{g}_{1(2r+1)}=\{x^u\partial_i||u|=2r, i\in\mathrm{I}_1\}$.

\subsection{Exceptional derivations}
Throughout this section assume that $|\omega|=n$ is odd.
In this section we shall consider three series of the so-called exceptional derivations from $H_{\overline{0}}$ into $W_{\overline{1}}$.
And we will see that these exceptional derivations are all outer.

Let us define the first series of exceptional derivations from $H_{\overline{0}}$ into $W_{\overline{1}}$.
For $i\in\mathrm{I}_0$ and $q\in\mathbb{N}$, define
$$\Phi_{i}^{(q)}:H_{\overline{0}}\rightarrow W_{\overline{1}},~~~~~\mathrm{D_H}(f)\mapsto\partial_i^{p^q}(f)\mathrm{D_H}(x^{\omega}).$$
Since $\mathrm{Ker(D_H)}=\mathbb{F}1$, note that $\Phi_{i}^{(q)}$ is well defined.

\blem\label{lem41}
~$\Phi_{i}^{(q)}\in\mathrm{Der}(H_{\overline{0}},W_{\overline{1}})$ and $\mathrm{ad}(\Phi_{i}^{(q)})=n-p^q$.
\elem

Now we define the second series of exceptional derivations.
We have known that $(\mathrm{ad}\partial_i)^{p^q}$ is a derivation of $H_{\overline{0}}$.
For $i\in\mathrm{I}_0$ and $q\in\mathbb{N}$, define
$$\Theta_{i}^{(q)}:H_{\overline{0}}\rightarrow W_{\overline{1}},~~~~~\mathrm{D_H}(f)\mapsto x^{\omega}(\mathrm{ad}\partial_i)^{p^q}(\mathrm{D_H}(f)).$$
Moreover, we have the following

\blem\label{lem42}
~$\Theta_{i}^{(q)}\in\mathrm{Der}(H_{\overline{0}},W_{\overline{1}})$ and $\mathrm{ad}(\Theta_{i}^{(q)})=n-p^q$.
\elem
The proofs of Lemmas \ref{lem41} and \ref{lem42} are analogous to \cite{Liu&Su&Zhang} Propositions 4.1 and 4.2.

Let us consider the third series of exceptional derivations.
For $i\in (1,\cdots,r)\cup(m+1,\cdots,m+s)$,
define
$$\Psi^{(i)}:H_{\overline{0}}\rightarrow W_{\overline{1}},~~~~~\mathrm{D_H}(f)\mapsto \partial_i\partial_{i'}(f)\mathrm{D_H}(x^{\omega})~~~~\mathrm{for}~f\in\mathcal{O}(m,n;\underline{t}).$$

\blem
~$\Psi^{(i)}\in\mathrm{D_H}(H_{\overline{0}},W_{\overline{1}})$ and $\mathrm{ad}(\Psi^{(i)})=n-2$.
\elem

\bpf
For $i\in (1,\cdots,r)$, the proof is analogous to \cite{Liu&Su&Zhang} Proposition 4.3.
We want to verify the conclusion for $i\in (m+1,\cdots,m+s)$.
For $\alpha$, $\beta\in\mathbb{A}$, $u$, $v\in\mathbb{B}^0$,
the following equation holds
\begin{align}\label{eq41}
&\Psi^{(i)}([\mathrm{D_H}(x^{(\alpha)}x^u),\mathrm{D_H}(x^{(\beta)}x^v)])\notag\\
=
&[\Psi^{(i)}(\mathrm{D_H}(x^{(\alpha)}x^u)),\mathrm{D_H}(x^{(\beta)}x^v)]+
[\mathrm{D_H}(x^{(\alpha)}x^u),\Psi^{(i)}(\mathrm{D_H}(x^{(\beta)}x^v))].
\end{align}

Case 1. $u\neq\emptyset$, $v\neq\emptyset$.
Since $|u|\geq2$, $|v|\geq2$ in this case, two sides of Eq.(\ref{eq41}) vanish.

Case 2. $u=\emptyset$, $v\neq\emptyset$.
The left-hand side of Eq.(\ref{eq41}) is as follows
\begin{align*}
\Psi^{(i)}([\mathrm{D_H}(x^{(\alpha)}),\mathrm{D_H}(x^{(\beta)}x^v)])=&\Psi^{(i)}(\mathrm{D_H}(\sum_{j=1}^{m}\tau(j)\partial_j(x^{(\alpha)})\partial_{j'}(x^{(\beta)})x^v))\\
=&\sum_{j=1}^{m}\tau(j)(x^{(\alpha-\varepsilon_j)})(x^{(\beta-\varepsilon_{j'})})\partial_i\partial_{i'}(x^v)\mathrm{D_H}(x^{\omega})
\end{align*}
The right-hand side of Eq.(\ref{eq41}) is as follows
\begin{align*}
&[\Psi^{(i)}(\mathrm{D_H}(x^{(\alpha)}x^u)),\mathrm{D_H}(x^{(\beta)}x^v)]+
[\mathrm{D_H}(x^{(\alpha)}x^u),\Psi^{(i)}(\mathrm{D_H}(x^{(\beta)}x^v))]\\
=&0+[\mathrm{D_H}(x^{(\alpha)}),x^{(\beta)}\partial_i\partial_{i'}(x^v)\mathrm{D_H}(x^{\omega})]\\
=&\sum_{j=1}^{m}\tau(j)(x^{(\alpha-\varepsilon_j)})(x^{(\beta-\varepsilon_{j'})})\partial_i\partial_{i'}(x^v)\mathrm{D_H}(x^{\omega})
\end{align*}
This proves Eq.(\ref{eq41}) in this case.

Case 3. $u\neq\emptyset$, $v=\emptyset$. The proof is analogous to Case 2.

Case 4. $u=v=\emptyset$. Obviously two sides of Eq.(\ref{eq41}) vanish.
\epf

\subsection{Derivations vanishing on the top}
Throughout this section assume that $|\omega|=n$ is odd unless otherwise stated.
Let us consider the elements in $\mathcal{M}$.

\blem\label{lem43}
Let $\phi\in\mathrm{Der}(\mathcal{J},W_{\overline{1}})$ be homogeneous such that $\phi(\mathcal{J}_{-1}\oplus\mathcal{J}_0)=0$.
Suppose that $\mathrm{zd}(\phi)$ is even and $\phi(\mathrm{D_H}(x^{(b\varepsilon_i)}))=0$ for all $b<a\leq\pi_i$,
where $a$ is a fixed positive integer and $i\in\mathrm{I}_0$.
Then the following statements hold:
\begin{itemize}
\item[\rm{(1)}] $a\equiv0$ $(\mathrm{mod}~p)$. If $a$ is not any $p-$power, then $\phi(\mathrm{D_H}(x^{(a\varepsilon_i)}))=0$. If $a=p^q$ for some $q\in\mathbb{N}$,~then there is $\lambda_i^{(q)}\in\mathbb{F}$ such that $$(\phi-\lambda_i^{(q)}\Phi_{i}^{(q)})(\mathrm{D_H}(x^{(a\varepsilon_i)}))=0.$$
\item[\rm{(2)}] $a\equiv1$ $(\mathrm{mod}~p)$. If $a-1$ is not any $p-$power, then $\phi(\mathrm{D_H}(x^{(a\varepsilon_i)}))=0$. If $a-1=p^q$ for some $q\in\mathbb{N}$,~then there is $\mu_i^{(q)}\in\mathbb{F}$ such that $$(\phi-\mu_i^{(q)}\Theta_{i}^{(q)})(\mathrm{D_H}(x^{(a\varepsilon_i)}))=0.$$
\item[\rm{(3)}]In the other conditions, $\phi(\mathrm{D_H}(x^{(a\varepsilon_i)}))=0$.
\end{itemize}
\elem

\bpf
$\rm{(1)}$ When $a\equiv0$ $(\mathrm{mod}~p)$, by Lemmas \ref{lem5} and \ref{lem7}~(1),
we may assume that
\begin{eqnarray}\label{eq42}
\phi(\mathrm{D_H}(x^{(a\varepsilon_i)}))=\sum_{r\in\mathrm{I_1}}                   \alpha_rx^{\widehat{r'}}\partial_r,~~~~~~~~~~~\mathrm{where}~\alpha_r\in\mathbb{F}.
\end{eqnarray}

If $a$ is not any $p-$power, then $a$ is written as the $p-$adic form $a=\sum_{r=1}^tc_rp^r$, $c_t\neq0$.
Obviously, ${a\choose p^{t}}\neq0$ ($\mathrm{mod}~p$).
For $i\in\mathrm{I}_0$, we have
\begin{eqnarray}\label{eq43}
[\mathrm{D_H}(x^{(p^t\varepsilon_i)}x_{i'}),\mathrm{D_H}(x^{((a-p^t+1)\varepsilon_i)})]=\tau(i'){a\choose p^{t}}\mathrm{D_H}(x^{(a\varepsilon_i)}).
\end{eqnarray}
Since $a-p^t+1<a$, there is the equation that $\phi(\mathrm{D_H}(x^{((a-p^t+1)\varepsilon_i)}))=0$.
On the other hand, we have
$$[\mathrm{D_H}(x^{(2\varepsilon_{i'})}),\mathrm{D_H}(x^{(a\varepsilon_i)})]=\tau(i')\mathrm{D_H}(x^{((a-1)\varepsilon_i)}x_{i'}).$$
It follow from Eq.(\ref{eq42}) that $\phi(\mathrm{D_H}(x^{((a-1)\varepsilon_i)}x_{i'}))=0$.
Since $p^t<a-1$, it follows that
$$\phi(\mathrm{D_H}(x^{(p^t\varepsilon_i)}x_{i'}))=\tau(i')[\mathrm{D_H}(x^{(2\varepsilon_{i'})},\phi(\mathrm{D_H}(x^{((p^t+1)\varepsilon_i)}))]=0. $$
Then $\phi(\mathrm{D_H}(x^{(p^t\varepsilon_i)}x_{i'}))=0$.
It is obtained that $\phi(\mathrm{D_H}(x^{(a\varepsilon_i)}))=0$.

If $a=p^q$ for some $q\in\mathbb{N}$, then $\mathrm{ad}(\phi)+a$ is odd.
Since
$[\mathrm{D_H}(x^kx^l),\mathrm{D_H}(x^{(a\varepsilon_i)})]=0$ for $k,l\in\mathrm{I}_1$,
applying $\phi$ on the equation, we have
$$(\alpha_{l'}x^l\partial_{k'}(x^{\widehat{l}})+\alpha_kx^{\widehat{k'}})\partial_{l'}-(\alpha_{k'}x^k\partial_{l'}(x^{\widehat{k}})+\alpha_lx^{\widehat{l'}})\partial_{k'}=0.$$
Furthermore,
$\alpha_{l'}x^l\partial_{k'}(x^{\widehat{l}})+\alpha_kx^{\widehat{k'}}=0$.
Without loss of generality, assume that $k<l'$.
It follows that  $\alpha_{k'}(-1)^{k'+l-3}+\alpha_{l}=0$,
that is, $\alpha_{l}=(-1)^{k'+l}\alpha_{k'}$.
Recall that $a=p^q$.
Put $\lambda_i^{(q)}:= c_{m+1}$.
Then we obtain from  Eq.(\ref{eq42}) that
$$\phi(\mathrm{D_H}(x^{(a\varepsilon_i)}))=\lambda_i^{(q)}\mathrm{D_H}(x^{\omega}).$$

$\rm{(2)}$
When $a\equiv1$ $(\mathrm{mod}~p)$, by Lemmas \ref{lem5} and \ref{lem7}~(2),
we may assume that
\begin{eqnarray}\label{eq44}
\phi(\mathrm{D_H}(x^{(a\varepsilon_i)}))=\mu_{i}x^{\overline{u}_i}x^{m+n}\partial_{i'},
\end{eqnarray}
where $\mu_{i}\in\mathbb{F}$ and $x^l,~x^{l'}$ are both in $x^{\overline{u}_i}$ for $l\in\{m+1,\cdots,m+s\}$.

If $a-1$ is not any $p-$power, then $a-1$ is written as the $p-$adic form $a-1=\sum_{r=1}^tc_rp^r$, $c_t\neq0$.
Obviously, ${a\choose p^{t}}\neq0$ (mod$p$).
Just as the proof of (1), we have $\phi(\mathrm{D_H}(x^{(a\varepsilon_i)}))=0$.

If $a-1=p^q$ for some $q\in\mathbb{N}$, then $\mathrm{ad}(\phi)+a$ is even.
For $k\in\mathrm{I}_1\setminus\{m+n\}$, we obtain that
$$0=[\mathrm{D_H}(x^kx^{m+n}),\mathrm{D_H}(x^{(a\varepsilon_i)})]=
-x^k\partial_{m+n}(x^{\overline{u}_i}x^{m+n})\partial_{i'}. $$
This implies that $x^{\overline{u}_i}x^{m+n}=x^{\omega}$.
Therefore, we can prove
$$\phi(\mathrm{D_H}(x^{(a\varepsilon_i)}))=\mu_ix^{\omega}\partial_{i'}.$$

$\rm{(3)}$ It is direct to prove by Lemmas \ref{lem5} and \ref{lem7}.
\epf

\blem\label{lem44}
Let $\phi\in\mathrm{Der}(\mathcal{J},W_{\overline{1}})$ be homogeneous such that $\phi(\mathcal{J}_{-1}\oplus\mathcal{J}_0)=0$.
Suppose that $\mathrm{zd}(\phi)$ is odd and $\phi(\mathrm{D_H}(x^{(b\varepsilon_i)}))=0$ for all $b<a\leq\pi_i$,
where $a$ is a fixed positive integer and $i\in\mathrm{I}_0$.
Then $\phi(\mathrm{D_H}(x^{(a\varepsilon_i)}))=0$.
\elem

\bpf
We also consider the three cases of $a$.
If $a\equiv0$ $(\mathrm{mod}p)$ and $a$ is not any $p-$power,
we can prove $\phi(\mathrm{D_H}(x^{(a\varepsilon_i)}))=0$ in the similar method of Lemma \ref{lem43} (1).
If $a=p^q$ for some $q\in\mathbb{N}$, then $\mathrm{zd}(\phi)+a$ is even.
So we can obtain $\phi(\mathrm{D_H}(x^{(a\varepsilon_i)}))\in\mathrm{E}(\mathfrak{g}_1)$,
which is a contradiction with Lemma \ref{lem7} (1).
Then $\phi(\mathrm{D_H}(x^{(a\varepsilon_i)}))=0$.
If $a\equiv1$ $(\mathrm{mod}p)$ and $a-1$ is not any $p-$power,
we can prove $\phi(\mathrm{D_H}(x^{(a\varepsilon_i)}))=0$ in the similar method of Lemma \ref{lem43} (2).
If $a-1=p^q$ for some $q\in\mathbb{N}$, then $\mathrm{zd}(\phi)+a$ is odd.
We can obtain $\phi(\mathrm{D_H}(x^{(a\varepsilon_i)}))\in\mathrm{O}(\mathfrak{g}_1)$,
which is a contradiction with Lemma \ref{lem7} (2).
Then $\phi(\mathrm{D_H}(x^{(a\varepsilon_i)}))=0$.
In the other conditions, it is direct to prove by Lemmas \ref{lem5} and \ref{lem7}.
\epf

Put $\mathfrak{D}=\sum_{t\in\mathrm{I}_1}x^tx^{m+n}\partial_t$.
For the elements in $\mathcal{N}$, we have the following lemma.
\blem
Let $\phi\in\mathrm{Der}(\mathcal{J},W_{\overline{1}})$ be homogeneous such that $\phi(\mathcal{J}_{-1}\oplus\mathcal{J}_0)=0$.
Then there is $\lambda\in \mathbb{F}$ such that $(\phi-\lambda\mathrm{ad}\mathfrak{D})(\mathcal{N})=0$.
\elem

\bpf
Recall that $\mathcal{N}=\{\mathrm{D_H}(x_ix^u)|i\in\mathrm{I}_0,~u\in\mathbb{B}_2\}$.
By Lemma \ref{lem7} (3),~we obtain that
$$\phi(\mathrm{D_H}(x_ix^kx^l))=\alpha_{i,\langle k,l\rangle}x^{k}x^{l}x^{m+n}\partial_{i'},$$
where $i\in\mathrm{I}_0$, $k,l\in\mathrm{I}_1$ and $\alpha_{i,\langle k,l\rangle}\in\mathbb{F}$.
Note that $\phi(\mathcal{J}_0)=0$.
For $s\in \mathrm{I}_1\setminus{\{l'\}}$, one can get the equation that $[\mathrm{D_H}(x^{k'}x^s),\mathrm{D_H}(x_ix^kx^l)]=\mathrm{D_H}(x_ix^sx^l)$.
Applying $\phi$ to the equation,
we obtain $\alpha_{i,\langle k,l\rangle}=\alpha_{i,\langle s,l\rangle}$.
This proves that $\alpha_{i,\langle k,l\rangle}$ is only dependent on the choice of $i$.
We denote $\alpha_{i}=\alpha_{i,\langle k,l\rangle}$ for $i\in\mathrm{I}_0$.
Then $\phi(\mathrm{D_H}(x_ix^u))=\alpha_{i}x^ux^{m+n}\partial_{i'}$ for all $x^u\in\mathbb{B}_{2}$.
For $i,~j\in\mathrm{I}_0$, we have $[\mathrm{D_H}(x_ix_{j'}),\mathrm{D_H}(x_jx^{u})]=\tau(j')(1+\delta_{i,j'})\mathrm{D_H}(x_ix^{u})$.
Applying $\phi$ to the equation we have $\tau(i)\alpha_{i}=\tau(j)\alpha_{j}$ for $i,~j\in\mathrm{I}_0$.
Let $\lambda=\frac{1}{2}\tau(i)\alpha_{i}$. Then
$$(\phi-\lambda\mathrm{ad}\mathfrak{D})(\mathrm{D_H}(x_ix^u))=\alpha_ix^u\partial_{i'}-2\tau(i)\lambda x^u\partial_{i'}=0.$$
The proof is complete.
\epf

\bprop\label{lem45}
Let $\phi\in\mathrm{Der}(\mathcal{J},W_{\overline{1}})$ be homogeneous and $\phi(\mathcal{J}_{-1}\oplus\mathcal{J}_0)=0$.
Then there are $\lambda$,~$\lambda_{r}^{(s_r)}$ and $\mu_{r}^{(s_r)}\in\mathbb{F}$, where $r\in\mathrm{I_0}$ and $1\leq s_r< t_r$, such that
$$\phi=\lambda(\mathrm{ad}\mathfrak{D})+\sum_{r\in\mathrm{I_0}}\sum_{1\leq s_r\leq t_r-1}(\lambda_{r}^{(s_r)}\Phi_{i}^{(s_r)}+\mu_{r}^{(s_r)}\Theta_{i}^{(s_r)}).$$
\eprop

\bpf In view of Lemmas 4.4 and 4.5, there are $\lambda_{r}^{(s_r)}$ and $\mu_{r}^{(s_r)}\in\mathbb{F}$, and we can prove by induction on $a_i$ that
$$\left(\phi-\sum_{r\in\mathrm{I_0}}\sum_{1\leq s_r\leq t_r-1}(\lambda_{r}^{(s_r)}\Phi_{i}^{(s_r)}+\mu_{r}^{(s_r)}\Theta_{i}^{(s_r)})\right)(\mathrm{D_H}(x^{(a_i\varepsilon_i)}))=0,$$
for all $i\in\mathrm{I}_0$, $1\leq a_i\leq\pi_i$.
By Lemma 4.6, there is $\lambda\in \mathbb{F}$ such that
$$\varphi=\phi-\sum_{r\in\mathrm{I_0}}\sum_{1\leq s_r\leq t_r-1}(\lambda_{r}^{(s_r)}\Phi_{i}^{(s_r)}+\mu_{r}^{(s_r)}\Theta_{i}^{(s_r)})-\lambda(\mathrm{ad}\mathfrak{D})$$
vanish on $\mathcal{N}$.
It is easy to see that $\varphi$ vanish on $\mathcal{M}$ and Lemma 3.2 ensures that $\varphi=0$.
The remaining is clear and the proof is complete.
\epf

\bthm
Let $\phi\in\mathrm{Der}(H_{\overline{0}},W_{\overline{1}})$ be homogeneous and $\phi((H_{\overline{0}})_{-1}\oplus(H_{\overline{0}})_0)=0$. If $n$ is even,
 then there are $\mu$, $\lambda$,~$\lambda_{r}^{(s_r)}$ and $\mu_{r}^{(s_r)}\in\mathbb{F}$, where $r\in\mathrm{I_0}$ and $1\leq s_r< t_r$, such that
$$\phi=\Gamma_{\mu}+\lambda(\mathrm{ad}\mathfrak{D})+\sum_{r\in\mathrm{I_0}}\sum_{1\leq s_r\leq t_r-1}(\lambda_{r}^{(s_r)}\Phi_{i}^{(s_r)}+\mu_{r}^{(s_r)}\Theta_{i}^{(s_r)}).$$
\ethm
\bpf
It is direct to obtain the consequence of Propositions 3.3 and 4.7 .
\epf
Obviously, we can get the following conclusion for the even integer $n$ by Proposition 3.3 and Lemma 3.7.
\bthm
 Let $\phi\in\mathrm{Der}(H_{\overline{0}},W_{\overline{1}})$ be homogeneous and $\phi((H_{\overline{0}})_{-1}\oplus(H_{\overline{0}})_0)=0$. If $n$ is even, then $\phi=0$.
\ethm

\subsection{Homogeneous derivation}

\blem$^{[13]}$
Suppose that $L$ is a $\mathbb{Z}$-graded subalgebra of~$W_{\overline{0}}$ satisfying $L_{-1}=(W_{\overline{0}})_{-1}$.
Let $\phi\in \mathrm{Der}(L,W_{\overline{1}})$ with $\mathrm{zd}(\phi)=t$.
Then there is $E\in(W_{\overline{1}})_t$ such that $$(\phi-\mathrm{ad}E)(L_{-1})=0.$$
\elem

\blem
Let $\phi\in\mathrm{Der}(\mathcal{J},W_{\overline{1}})$ be $\mathbb{Z}$-homogeneous with the odd degree such that $\phi(\mathcal{J}_{-1})=0$.
Then there is $\lambda_i\in\mathbb{F}$ such that
$(\phi-\lambda_i\Psi^{(i)})(\mathcal{J}_0)=0$
\elem

\bpf
Since $\mathrm{zd}(\phi)$ is odd, $\mathrm{zd}(\phi(\mathcal{J}_0))$ is odd.
We can obtain $\phi(\mathcal{J}_0)\in\mathrm{O}(\mathfrak{g}_1)$.
For $i\in\mathrm{I_0}$, by Lemmas \ref{lem5} and \ref{lem7}~(4), we can assume that
\begin{eqnarray}\label{eq100}
\phi(\mathrm{D_H}(x_ix_{i'}))=\sum_{r\in\mathrm{I_1}}                   \alpha_rx^{\widehat{r'}}\partial_r,~~~~~~~~~~~\mathrm{where}~\alpha_r\in\mathbb{F}.
\end{eqnarray}
For arbitrary $k,~l\in\mathrm{I_1}$, with $k\neq l'$, we have
$$0=[\mathrm{D_H}(x_kx_l),\phi(\mathrm{D_H}(x_ix_{i'}))]=\sum_{r\in\mathrm{I_1}}\alpha_r(x_l\partial_{k'}-x_k\partial_{l'})x^{\widehat{r'}}\partial_r-
\alpha_lx^{\widehat{l'}}\partial_{k'}+\alpha_rx^{\widehat{k'}}\partial_{l'}.$$
Furthermore,
$\alpha_{l'}x^l\partial_{k'}(x^{\widehat{l}})+\alpha_kx^{\widehat{k'}}=0$.
Without loss of generality, assume that $k<l'$.
It follows that  $\alpha_{k'}(-1)^{k'+l-3}+\alpha_{l}=0$,
that is, $\alpha_{l}=(-1)^{k'+l}\alpha_{k'}$.
Put $\lambda_i:= c_{m+1}$.
Then we obtain from  Eq.(\ref{eq100}) that
$$\phi(\mathrm{D_H}(x_ix_{i'}))=\lambda_i\mathrm{D_H}(x^{\omega}).$$

For $i,j\in\mathrm{I_0}$, with $j\neq i'$, by Lemmas \ref{lem5} and \ref{lem7}~(5),  we can assume that
$$\phi(\mathrm{D_H}(x_ix_j))= \alpha_{i'}x^{\overline{u}_{i'}}x^{m+n}\partial_{i'}+\alpha_{j'}x^{\overline{u}_{j'}}x^{m+n}\partial_{j'}.$$
Then there is a contradiction. So $\phi(\mathrm{D_H}(x_ix_j))=0$.

By Lemmas \ref{lem5}, \ref{lem7}~(6) and (7), for $k.~l\in\mathrm{I_1}$ with $l\neq k'$, we can assume that
\begin{align}\label{eq200}
&\phi(\mathrm{D_H}(x_kx_{k'}))=\sum_{r\in\mathrm{I_1}}\alpha_rx^{\widehat{r'}}\partial_r,
\\
&\phi(\mathrm{D_H}(x_kx_{l}))=\beta_{k'}x^{\overline{u}_{k'}}\partial_{k'}+\beta_{l'}x^{\overline{u}_{l'}}\partial_{l'},
\end{align}
where $\alpha_r,~\beta_{k'},~\beta_{l'}\in\mathbb{F}$ and $x^l,~x^{l'}$ are both in $x^{\overline{u}_{k'}}$ and $x^{\overline{u}_{l'}}$ for $l\in\{m+1,\cdots,m+s\}$.
Since
$[\mathrm{D_H}(x_kx_{k'}),\mathrm{D_H}(x_sx_t)]=0$ for $s,~t\in\mathrm{I}_1\setminus\{k,~k'\}$,
applying $\phi$ on the equation, we have that
$$0=[\sum_{r\in\mathrm{I_1}}\alpha_rx^{\widehat{r'}}\partial_r,x_t\partial_{s'}-x_s\partial_{t'}]
=\alpha_tx^{\widehat{t'}}\partial_{s'}-\alpha_sx^{\widehat{s'}}\partial_{t'}
-\alpha_{t'}x_t\partial_{s'}(x^{\widehat{t}})\partial_{t'}+\alpha_{s'}x_s\partial_{t'}(x^{\widehat{s}})\partial_{s'}.
$$
By the calculation of the coefficients of $\partial_{t'}$ and $\partial_{s'}$,
we have $\alpha_t=(-1)^{t+s}\alpha_s$ for $s,~t\in\mathrm{I}_1\setminus\{k,~k'\}$.
We will particularly discuss $\alpha_k$ and $\alpha_{k'}$.
For $t\in\mathrm{I_1}\setminus\{k'\}$, we have the equation that
$$[\mathrm{D_H}(x_kx_{k'}),\mathrm{D_H}(x_kx_t)]=-\mathrm{D_H}(x_kx_t).$$
Applying $\phi$ on the equation, by Eq.(4.7),
we have
$$\alpha_tx^{\widehat{t'}}\partial_{k'}-\alpha_kx^{\widehat{k'}}\partial_{t'}+\alpha_{k'}x_k\partial_{t'}(x^{\widehat{k}})\partial_{k'} +\alpha_{t'}x_t\partial_{k'}(x^{\widehat{t}})\partial_{t'}=\beta_{k'}x^{\overline{u}_{k'}}\partial_{k'}+\beta_{t'}x^{\overline{u}_{t'}}\partial_{t'}.
$$
Since $x^l,~x^{l'}$ are both in $x^{\overline{u}_{k'}}$ and $x^{\overline{u}_{l'}}$ for $l\in\{m+1,\cdots,m+s\}$, we have that $\alpha_k=(-1)^{t'+k}\alpha_{t'}$ and  $\alpha_{k'}=(-1)^{t+k'}\alpha_{t}$.
Furthermore, we also get $\beta_{k'}=\beta_{t'}=0$,
for $k, t\in\mathrm{I}_1$, with $t\neq k'$.
Then we obtain that $\phi(\mathrm{D_H}(x_kx_t))=0$, for $k, t\in\mathrm{I}_1$, with $t\neq k'$.
Put $\lambda_k:= c_{m+1}$.
Then we obtain from Eq.(4.6) that
$$\phi(\mathrm{D_H}(x_kx_{k'}))=\lambda_k\mathrm{D_H}(x^{\omega}).$$
For $i\in(1,\cdots,r)\cup(m+1,\cdots,m+s)$, we put $$\varphi=\phi-\lambda_i\Psi^{(i)},$$
where $\lambda_i\in\mathbb{F}$.
By the above proof, obviously, $\varphi(\mathcal{J}_0)=0$.
\epf

\blem
Let $\phi\in\mathrm{Der}(\mathcal{J},W_{\overline{1}})$ be $\mathbb{Z}$-homogeneous with the even degree such that $\phi(\mathcal{J}_{-1})=0$.
Then $\phi(\mathcal{J}_0)=0$.
\elem

\bpf
Since $\mathrm{zd}(\phi)$ is even, $\mathrm{zd}(\phi(\mathcal{J}_0))$ is even.
We can obtain $\phi(\mathcal{J}_0)\in\mathrm{E}(\mathfrak{g}_1)$.
By Lemmas \ref{lem5} and \ref{lem7}, for $i\in\mathrm{I}_0,~k,~l\in\mathrm{I}_1$ we can obtain
that $\phi$ vanish on $\mathrm{D_H}(x_ix_{i'})$, $\mathrm{D_H}(x_kx_{k'})$ and $\mathrm{D_H}(x_kx_{l})$.
For $i,j\in\mathrm{I}_0$, with $j\neq i, i'$, by Lemma \ref{lem7}~(5), we can assume that
\begin{align}\label{eq215}
\phi(\mathrm{D_H}(x_ix_j))=\alpha_{i'}x^{\overline{u}_{i'}}x^{m+n}\partial_{i'}+\alpha_{j'}x^{\overline{u}_{j'}}x^{m+n}\partial_{j'}.
\end{align}
Since $[\mathrm{D_H}(x_ix_j),\mathrm{D_H}(x_{i'}x_{j'})]=\tau(i)\mathrm{D_H}(x_jx_{j'})+\tau(j)\mathrm{D_H}(x_ix_{i'})$,
applying $\phi$ on the equation,
we can obtain that $\alpha_{i}=\alpha_{j}=\alpha_{i'}=\alpha_{j'}=0$ by Eq.(4.8).
Then the conclusion is proved.
\epf


\bprop\label{lem405}
Let $H_{\overline{0}}$ be the even parts of $H$ and $W_{\overline{1}}$ the odd parts of $W$.
$\mathcal{J}$~is an idea of $H_{\overline{0}}$ of codimension 1.
Then the following conclusion holds:
\begin{displaymath}
\mathrm{Der}(\mathcal{J},W_{\overline{1}})= \left\{ \begin{array}{ll}
\mathrm{ad}W_{\overline{1}}+\mathbb{F}\Psi^{(i)} +\sum_{r\in\mathrm{I_0}}\sum_{1\leq s_r\leq t_r-1}(\mathbb{F}\Phi_{i}^{(s_r)}+\mathbb{F}\Theta_{i}^{(s_r)})   & \textrm{$n$ is odd,}\\
\mathrm{ad}W_{\overline{1}} & \textrm{$n$ is even.}
\end{array} \right.
\end{displaymath}
%
%
\eprop

\bpf
It is direct to prove by Proposition 4.7, Lemmas 4.9, 4.10 and 4.11.
\epf

Obviously, we can get the following conclusion by Theorems 4.8, 4.9 and Proposition 3.3.

\bthm
Let $H_{\overline{0}}$ be the even parts of $H$ and $W_{\overline{1}}$ the odd parts of $W$.
Then the following conclusion holds:
\begin{displaymath}
\mathrm{Der}(H_{\overline{0}},W_{\overline{1}})= \left\{ \begin{array}{ll}
\mathbb{F}\Gamma_{1}+\mathrm{ad}W_{\overline{1}}+\mathbb{F}\Psi^{(i)} +\sum_{r\in\mathrm{I_0}}\sum_{1\leq s_r\leq t_r-1}(\mathbb{F}\Phi_{i}^{(s_r)}+\mathbb{F}\Theta_{i}^{(s_r)})   & \textrm{$n$ is odd,}\\
\mathrm{ad}W_{\overline{1}} & \textrm{$n$ is even.}
\end{array} \right.
\end{displaymath}
\ethm

\end{document}